# Поиск стохастических равновесий в транспортных моделях равновесного распределения потоков

*Гасников А.В., Гасникова Е.В.,
Двуреченский П.Е., Ершов Е.И., Лагуновская А.А.*

(Центр исследований транспортной политики, Институт экономики транспорта и транспортной политики, НИУ ВШЭ; ИППИ РАН; ПреМоЛаб МФТИ)

**Аннотация**

В работе предложены эффективные способы поиска стохастических равновесий в популяционных играх загрузок. Поиск равновесия Нэша в таких играх всегда сводится к задаче оптимизации. Мы рассматриваем модели равновесного распределения потоков по путям Бэкмана и Нестерова–де Пальмы. Поиск стохастических равновесий Нэша(–Вардропа) приводит к энтропийной регуляризации выпуклого функционала, отвечающего этим моделям. Данная работа посвящена тому, как эффективно решать такого рода задачи. В основе подхода лежит идея композитной оптимизации и особенность постановки, что функционал имеет вид суммы (сепарабельный функционал). Это обстоятельство, вместе с неограниченностью константы Липшица градиента функционала, мотивирует переформулировку исходной задачи оптимизации таким образом, чтобы этот сепарабельный функционал стал композитным членом. Рассматриваются и развиваются также и классические способы решения отмеченной задачи с помощью аппарата характеристических функций на графе.

**Ключевые слова:** Композитный быстрый градиентный метод, разреженность, равновесное распределение потоков.

## 1. Введение

В работе [1] было анонсировано, что в цикле статей планируется описать эффективные численные методы поиска равновесий (и стохастических равновесий) в популярных моделях распределения транспортных потоков по путям. В частности, это планировалось сделать для модели Бэкмана и модели стабильной динамики (Нестерова–деПальмы). В работе [2] мы сконцентрировались на поиске обычных равновесий (Нэша). Данная работа посвящена эффективным численным методам поиска стохастических равновесий в этих моделях. Стохастические равновесия играют важную роль в практических приложениях этих моделей [3, 4]. В популяционной теории игр, исследуемые нами стохастические равновесия иногда называют логит равновесиями [5]. Такие равновесия получаются, если вместо динамики наилучших ответов (приводящей к равновесию Нэша при условиях потенциальности игры, с выпуклым потенциалом [5]) и положительно коррелирующих с ней динамик (также приводящих к равновесию Нэша) рассматривается логит динамика (гиббсовская динамика), отражающая поведение ограниченно рациональных агентов, т.е. агентов действующих в условиях не полной информации (немного подробнее об этом будет написано в п. 2, посвященном постановке задачи и в п. 5). Последнее обстоятельство позволяет связать такие равновесия с достаточно популярным разделом современной экономической теории – с теорией дискретного выбора [6].

В пп. 2 – 4 приводятся новые подходы к поиску стохастических равновесий, в основе которых лежат современные результаты из композитной оптимизации [7, 8]. Ожидается, что эти методы будут неплохо работать в случае, когда число возможных путей в графе транспортной сети не очень большое по сравнению с числом ребер. В общем же случае число возможных путей может быть "экспоненциально" велико по сравнению с числом ребер. Тогда методы из пп. 2 – 4 доминируются методом, базирующимся на аппарате характеристических функций на графе [4, 9]. Дополнительно к описанию подхода работ [4, 9] мы приводим в п. 5 способ ускорения описанного там метода для графов специальной (упорядоченной) структуры, а также приводим новый рандомизированный метод, тесно связанный с алгоритмом из п. 3 [2].

В данной работе мы сосредоточились на двух моделях равновесного распределения транспортных потоков по путям: модели Бэкмана и модели стабильной динамики. Однако легко перенести все, что далее будет написано, например, на смешанные модели [2]. Также заложенные в этой статье и статье [2] подходы впоследствии должны позволить разработать эффективные численные методы для поиска равновесий в многостадийных транспортных моделях, одним из блоков которых является изучаемые в данной статье модели распределения транспортных потоков по путям [1, 10–12].

**2. Постановка задачи**

Рассмотрим транспортную сеть, которую будем представлять ориентированным графом $\langle V, E \rangle$, где $V$ – множество вершин (как правило, можно считать, что $|E|/4 \le |V| \le |E|$), а $E$ – множество ребер ($|E| = m$). Обозначим множество пар $w = (i, j)$ источник-сток через $OD$, $d_w$ – корреспонденция, отвечающая паре $w$, $x_p$ – поток по пути $p$; $P_w$ – множество путей, отвечающих корреспонденции $w$ (начинающихся в $i$ и заканчивающихся в $j$), $P = \bigcup_{w \in OD} P_w$ – множество всех путей ($|P| = n$). Затраты на прохождения ребра $e \in E$ описываются функцией $\tau_e(f_e)$, где $f_e$ – поток по ребру $e$.

Опишем марковскую логит динамику (также говорят гиббсовскую динамику) в повторяющейся игре загрузки графа транспортной сети [5]. Пусть каждой корреспонденции отвечает $d_w M$ агентов ($M \gg 1$), $\tau_e(f_e) := \tau_e(f_e/M)$. Пусть имеется $TN$ шагов ($N \gg 1$). Каждый агент независимо от остальных на шаге $t+1$ выбирает с вероятностью

$$\frac{\lambda}{N} \frac{\exp(-G_p^t/\gamma)}{\sum_{\tilde{p} \in P_w} \exp(-G_{\tilde{p}}^t/\gamma)}.$$

путь $p \in P_w$, где $G_p^t$ – затраты на пути $p$ на шаге $t$ ($G_p^0 \equiv 0$), а с вероятность $1 - \lambda/N$ путь, который использовал на шаге $t$. Такая динамика отражает ограниченную рациональность агентов (см. п. 5), и часто используется в популяционной теории игр [5] и теории дискретного выбора [6]. Оказывается эта марковская динамика в пределе $N \to \infty$ превращается в марковскую динамику в непрерывном времени (вырождающуюся при $\gamma \to 0+$ в динамику наилучших ответов [5]), которая в свою очередь допускает два предельных перехода (обоснование перестановочности этих пределов см. в [13]): $T \to \infty$, $M \to \infty$ или $M \to \infty$, $T \to \infty$. При первом порядке переходов мы сначала ($T \to \infty$) согласно эргодической теореме для марковских процессов (в нашем случае марковский процесс – модель стохастической химической кинетики с унарными реакциями в условиях детального баланса [14–16]) приходим к финальной (=стационарной) вероятностной мере, имеющей в основе мультиномиальное распределение. С ростом числа агентов ($M \to \infty$) эта мера концентрируется около наиболее вероятного состояния, поиск которого сводится к решению задачи (1) ниже. Функционал в этой задаче оптимизации с точностью до потенцирования и мультипликативных и аддитивных констант соответствует исследуемой стационарной мере – то есть это функционал Санова [16, 17]. При обратном порядке переходов, мы сначала осуществляем, так называемый, канонический скейлинг [13, 14], приводящий к детерминированной кинетической динамике, описываемой СОДУ на $x$, а затем ($T \to \infty$) ищем аттрактор получившейся СОДУ. Глобальным аттрактором оказывается неподвижная точка, которая определяется решением задачи (1) ниже. Более того, функционал, стоящий в (1), является функцией Ляпунова полученной кинетической динамики (то есть функционалом Больцмана). Последнее утверждение – достаточно общий факт (функционал Санова, является функционалом Больцмана), верный при намного более общих условиях [16].

Итак, рассматривается следующая задача поиска стохастического равновесия Нэша–Вардропа в модели Бэкмана равновесного распределения транспортных потоков по путям (см., например, [1, 3]):

$$\sum_{e \in E} \sigma_e(f_e) + \gamma \sum_{w \in OD} \sum_{p \in P_w} x_p \ln(x_p/d_w) \to \min_{f = \Theta x, x \in X}, \qquad (1)$$

где $\gamma > 0$; $\sigma_e(f_e) = \int_0^{f_e} \tau_e(z) dz$ – выпуклые функции;

$$\Theta = \|\delta_{ep}\|_{e \in E, p \in P} = \|\Theta^{\langle p \rangle}\|_{p \in P} = \|\Theta_e\|_{e \in E}^T, \quad \delta_{ep} = \begin{cases} 1, & e \in p; \\ 0, & e \notin p; \end{cases}$$

важно заметить, что матрица $\Theta$ сильно разрежена по столбцам ($s$ – среднее число ребер в пути, т.е. среднее число ненулевых элементов в столбце матрицы $\Theta$, как правило, можно считать, что $s \ll \sqrt{m}$);

$$X = \left\{ x \geq 0 : \sum_{p \in P_w} x_p = d_w, w \in OD \right\} \text{ – прямое произведение симплексов;}$$

в качестве $\tau_e(f_e)$ обычно выбирают BPR-функции [1, 4, 18, 19]:

$$\tau_e(f_e) = \bar{t}_e \cdot \left(1 + \rho \cdot (f_e/\bar{f}_e)^4\right).$$

В пределе модели стабильной динамики [1]

$$\tau_e^\mu(f_e) = \bar{t}_e \cdot \left(1 + \rho \cdot (f_e/\bar{f}_e)^{1/\mu}\right) \xrightarrow[\mu \to 0+]{} \begin{cases} \bar{t}_e, & 0 \leq f_e < \bar{f}_e \\ [\bar{t}_e, \infty), & f_e = \bar{f}_e \end{cases},$$

$$d\tau_e^\mu(f_e)/df_e \xrightarrow[\mu \to 0+]{} 0, \ 0 \leq f_e < \bar{f}_e$$

задача перепишется как [1]

$$\sum_{e \in E} f_e \overline{t}_e + \gamma \sum_{w \in OD} \sum_{p \in P_w} x_p \ln\left(x_p / d_w\right) \to \min_{\substack{f = \Theta x, \, x \in X \\ f \leq \overline{f}}}. \tag{2}$$

Эту задачу можно переписать как задачу энтропийно-линейного программирования, и использовать соответствующие численные методы [20, 21]. Это приводит к следующей оценке числа арифметических операций: $\mathrm{O}\left(sn\sqrt{L_2 R_2^2 / \varepsilon}\right)$, где $L_2 = \max_{p \in P} \left\|\Theta^{\langle p \rangle}\right\|_2^2 = H$ ($\Theta^{\langle p \rangle}$ – $p$-й столбец матрицы $\Theta$, $H$ – максимальное число ребер в пути, как правило, можно считать, что $H = \mathrm{O}\left(\sqrt{m}\right)$), $R_2$ – евклидов размер решения двойственной задачи (итерационный процесс в двойственном пространстве стартует с нуля), $\varepsilon$ – точность по функции (в прямой задаче). При этом ограничения выполняются с точностью $\varepsilon/R_2$. Интересно было бы попробовать использовать для решения двойственной задачи к (2) не быстрый градиентный метод FGM, как в [21], а его покомпонентный вариант ACRCD, также обладающий прямо-двойственной структурой [22]. Тем не менее, на данный момент нам не удалось теоретически обосновать, что это может принести какие-то дивиденды по сравнению с использованием обычного FGM (даже при дополнительных предположениях типа разреженности матрицы $\Theta$).

### 3. Композитный подход для задачи (1)

Прежде всего, заметим, что $f(x) = \Theta x$ и функции $\sigma_e(f_e)$ гладкие (с ограниченной константой Липшица градиента) на компакте $\{f = \Theta x : x \in X\}$. Следовательно, к задаче (1), записанной в следующей форме

$$\sum_{e \in E} \sigma_e \left(\Theta_e^T x\right) + \gamma \sum_{w \in OD} \sum_{p \in P_w} x_p \ln\left(x_p / d_w\right) \to \min_{x \in X}, \tag{1'}$$

можно применять композитный быстрый градиентный метод [7, 8], считая композитом энтропию $\gamma \sum_{w \in OD} \sum_{p \in P_w} x_p \ln\left(x_p / d_w\right)$, выбирая норму

$$\|x\| = \sqrt{\sum_{w \in OD} \left\|\{x_p\}_{p \in P_w}\right\|_1^2}$$

в прямом пространстве и 1-сильно выпуклую в этой норме прокс-функцию

$$d_1(x) = \sum_{w \in OD} d_w \cdot d^w\left(\{x_p\}_{p \in P_w}\right), \quad d^w\left(\{x_p\}_{p \in P_w}\right) = d_w \ln |P_w| + \sum_{p \in P_w} x_p \ln\left(x_p / d_w\right),$$

если $\gamma > 0$ мало (не учитываем сильную выпуклость функционала), иначе (учитываем сильную выпуклость функционала) выбираем 1-сильно выпуклую прокс-функцию

$$d_2(x) = \sum_{w \in OD} d_w \cdot \frac{1}{2(a_w - 1)} \left\|\{x_p\}_{p \in P_w}\right\|_{a_w}^2 \quad \text{с} \quad a_w = \frac{2 \ln |P_w|}{2 \ln |P_w| - 1}.$$

В случае, когда $\gamma > 0$ мало, имеем (согласно [7]) следующую оценку трудоемкость метода (общее число арифметических операций) будет (здесь использовано, что $\mathrm{O}(sN)$ – стоимость одной итерации, поскольку шаг композитного метода осуществляется по явным формулам, и основная сложность заключается в вычислении $f = \Theta x$)

$$\mathrm{O}\left(sn\sqrt{\frac{LR^2}{\varepsilon}}\right), \quad (3)$$

где

$$L = \max_{f=\Theta x,\, x\in X}\max_{\|h\|\le 1}\left\langle h, \Theta^T \mathrm{diag}\left\|\left\{\tau'_e(f_e)\right\}_{e\in E}\right\|\Theta h\right\rangle \le \max_{\substack{f=\Theta x,\, x\in X\\ e\in E}}\tau'_e(f_e)\cdot\left\|\left\{\max_{p\in P_w}\left\|\Theta^{\langle p\rangle}\right\|_2\right\}_{w\in OD}\right\|_2^2,$$

$$R^2 = \max_{x\in X} d_1(x) = \sum_{w\in OD} d_w^2 \ln|P_w|.$$

В случае, когда $\gamma > 0$ немало, имеем (согласно [7]) следующую оценку трудоемкость метода (здесь использовано, что $\mathrm{O}(sn)$ – стоимость одной итерации [21])

$$\mathrm{O}\left(sn\sqrt{\frac{L}{\mu}\chi}\cdot\ln\left(\frac{\mu\chi\|x_0-x_*\|^2}{\varepsilon}\right)\right), \quad (4)$$

где константа сильной выпуклости композита в введенной норме на $X$: $\mu = \gamma\cdot\left(\max_{w\in OD} d_w\right)^{-1}$, а множитель $\chi \le 2\max_{w\in OD}\ln|P_w|$ возникает из-за использования техники рестартов в прямом пространстве и прокс-функции $d_2(x)$ (см., например, стр. 56 [23][1] и п. 5 работы [24]). Из формул (3), (4) можно сделать вывод, что учитывать 1-сильную выпуклость композита на $X$ стоит, если $R^2/\varepsilon \gg \chi/\mu$, т.е. при

$$\gamma \gg \varepsilon\cdot\left(\sum_{w\in OD} d_w\right)^{-1}.$$

В частном случае $|OD|=1$ описанный здесь подход ранее предлагался в работе [21]. Можно использовать и другие способы решения задачи (1') из [21], например, базирующиеся на рандомизированном двойственном покомпонентном методе.

К сожалению, во всех подходах вылезает одно и то же слабое место: константа $L$ (даже при аккуратном оценивании) оказывается слишком большой. Связано это с поведением функций $\tau_e(f_e)$, имеющих быстро растущие производные. К сожалению, константа $L$ явно входит в размер шага метода и в критерий останова.

Помочь в решении отмеченной проблемы может следующее наблюдение. В действительности, функции $\tau_e(f_e)$ имеют штрафную природу: при приближении $f_e$ к пропускной способности $\bar{f}_e$ наблюдается резкий рост. На самом решении (в равновесии) и в его окрестностях не стоит ожидать, что пропускные способности сильно превышены. Это позволяет надеяться, что по ходу итерационного процесса это свойство также будет выполнено. Для этого нужно еще стартовать с разумной точки, в которой нет сильных нарушений пропускных способностей. По условию мы считаем, что это возможно, иначе искомое равновесие будет "одной большой пробкой". Если предположить, что есть метод, который настраивается на реальную константу Липшица градиента, а не на ее оценку сверху $L$, то во многом отмеченная проблема будет решена. Оказывается, можно предложить универсальный вариант [8, 25] описанного выше композитного метода (в том числе с учетом сильной выпуклости и возникающей неточности при вычислении шага градиентного отображения в сильно выпуклом случае [26]), который оптимально самонастраивается на текущую гладкость оптимизируемой функции (за вычетом композита).

---

[1] Заметим, что в сильно выпуклом случае мы используем сильную выпуклость композита (этим объясняется некоторое противоречие с тем, что написано в п. 3 [23] – там речь идет об одной и той же функции, когда говорится, что не имеет смысл рассматривать в сильно выпуклом случае прокс-структуры отличные от евклидовой: это верно для не композитных постановок).

Другой способ борьбы с большим $L$ будет описан далее.

### 4. Штрафное раздутие задач

Задачи (1) и (2) можно переписать следующим образом (см. [21] в связи с тем как следует выбирать $\lambda$)

$$\frac{1}{2}\|\Theta x - f\|_2^2 + \lambda_1 \cdot \left( \sum_{e \in E} \sigma_e(f_e) + \gamma \sum_{w \in OD} \sum_{p \in P_w} x_p \ln(x_p/d_w) \right) \to \min_{x \in X}, \qquad (1'')$$

$$\frac{1}{2}\|\Theta x - f\|_2^2 + \lambda_2 \cdot \left( \sum_{e \in E} f_e \bar{t}_e + \gamma \sum_{w \in OD} \sum_{p \in P_w} x_p \ln(x_p/d_w) \right) \to \min_{x \in X, f \leq \bar{f}}. \qquad (2'')$$

Мы выбираем норму в прямом пространстве $(x, f)$ следующим образом:

$$\|(x, f)\| = \|x\| + \|f\|_2.$$

Поскольку $\|(x, f)\|^2 \leq 2\|x\|^2 + 2\|f\|_2^2$, то прокс-функция выбирается $2d_1(x) + \|f\|_2^2$. Далее мы имеем композитную постановку с композитом вида

$$\sum_{e \in E} \sigma_e(f_e) + \gamma \sum_{w \in OD} \sum_{p \in P_w} x_p \ln(x_p/d_w).$$

С точностью до симплексных ограничений на $x$ это сепарабельный композит. Далее к нему следует применить (должным образом модернезировав) методы из работы [21].

Приведем, например, оценку трудоемкости композитного варианта быстрого градиентного метода для задач (1''), (2'')

$$\mathrm{O}\left( sn\sqrt{\frac{\breve{L}\breve{R}^2}{\varepsilon^2}} \right),$$

$$\breve{L} = \max_{\|h\| \leq 1} \langle h, \Theta^T \Theta h \rangle + \max_{\|f\|_2 \leq 1} \langle f, f \rangle = \left\| \left\{ \max_{p \in P_w} \|\Theta^{\langle p \rangle}\|_2 \right\}_{w \in OD} \right\|_2^2 + 1,$$

$$\breve{R}^2 = 2 \sum_{w \in OD} d_w^2 \ln|P_w| + \|\bar{f}\|_2^2,$$

обеспечивающую $\varepsilon$-решение (по функции) задачам (1), (2) и выполнение неравенства $\|\Theta x - f\|_2 \leq \varepsilon$. Отметим, что каждый шаг композитного метода считается по явным формулам, поскольку переменные распадаются на две группы: по $x$ это стандартное экспоненциальное взвешивание, а по $f$ потребуется решать уравнение четвертой степени для BPR-функций (см. п. 2), что также можно сделать по явным формулам.

### 5. Использование аппарата характеристических функций на графе

Запишем двойственную задачу к (1) [1, 2, 4] (далее мы используем обозначение $\operatorname{dom} \sigma^*$ – область определения сопряженной к $\sigma$ функции)

$$\min_{f,x}\left\{\sum_{e\in E}\sigma_e(f_e)+\gamma\sum_{w\in OD}\sum_{p\in P_w}x_p\ln(x_p/d_w):\ f=\Theta x,\ x\in X\right\}=$$

$$=\min_{f,x}\left\{\sum_{e\in E}\max_{t_e\in\operatorname{dom}\sigma_e^*}\left[f_e t_e-\sigma_e^*(t_e)\right]+\gamma\sum_{w\in OD}\sum_{p\in P_w}x_p\ln(x_p/d_w):\ f=\Theta x,\ x\in X\right\}=$$

$$=\max_{t\in\operatorname{dom}\sigma^*}\left\{\min_{f,x}\left[\sum_{e\in E}f_e t_e+\gamma\sum_{w\in OD}\sum_{p\in P_w}x_p\ln(x_p/d_w):\ f=\Theta x,\ x\in X\right]-\sum_{e\in E}\sigma_e^*(t_e)\right\}=$$

$$=-\min_{t\in\operatorname{dom}\sigma^*}\left\{\gamma\psi(t/\gamma)+\sum_{e\in E}\sigma_e^*(t_e)\right\}, \qquad (5)$$

где[2]

$$\psi(t)=\sum_{w\in OD}d_w\psi_w(t),\ \psi_w(t)=\ln\left(\sum_{p\in P_w}\exp\left(-\sum_{e\in E}\delta_{ep}t_e\right)\right),$$

$$f=-\nabla\psi(t/\gamma),\ x_p=d_w\frac{\exp\left(-\frac{1}{\gamma}\sum_{e\in E}\delta_{ep}t_e\right)}{\sum_{\tilde p\in P_w}\exp\left(-\frac{1}{\gamma}\sum_{e\in E}\delta_{e\tilde p}t_e\right)},\ p\in P_w,$$

для

$$\tau_e(f_e)=\bar t_e\cdot\left(1+\rho\cdot\left(\frac{f_e}{\bar f_e}\right)^{\frac{1}{\mu}}\right),$$

имеем [1] (BPR-функция, см. п. 2, получается при $\mu=1/4$)

$$\sigma_e^*(t_e)=\sup_{f_e\ge 0}\left((t_e-\bar t_e)\cdot f_e-\bar t_e\cdot\frac{\mu}{1+\mu}\cdot\rho\cdot\frac{f_e^{1+\frac{1}{\mu}}}{\bar f_e^{\frac{1}{\mu}}}\right)=\bar f_e\cdot\left(\frac{t_e-\bar t_e}{\bar t_e\cdot\rho}\right)^\mu\frac{(t_e-\bar t_e)}{1+\mu}.$$

---

[2] Обратим внимание, что вектор распределения потоков по путям $x$ при поиске стохастического равновесия ($\gamma>0$) получается не разреженным в отличие от поиска обычного равновесия Нэша(–Вардропа) [2] ($\gamma=0$). Как следствие, чтобы вычислить этот вектор требуются затраты существенно зависящие от потенциально огромной размерности $n$. К счастью, в приложениях, как правило, не требуется знание этого вектора, достаточно определить вектор потоков на ребрах $f$, который, как мы увидим ниже, может быть вычислен намного эффективнее (в частности, с затратами независящими от $n$).

Собственно, формула (5) есть не что иное, как отражение формулы $f = -\nabla \psi(t/\gamma)$ и связи $t_e = \tau_e(f_e)$, $e \in E$. Действительно, по формуле Демьянова–Данскина–Рубинова [1]

$$\frac{d\sigma_e^*(t_e)}{dt_e} = \frac{d}{dt_e} \max_{f_e \geq 0} \left\{ t_e f_e - \int_0^{f_e} \tau_e(z)dz \right\} = f_e : t_e = \tau_e(f_e).$$

В свою очередь, формула $f = -\nabla \psi(t/\gamma)$ может интерпретироваться, как следствие соотношений $f = \Theta x$ и формулы распределения Гиббса (логит распределения)

$$x_p = d_w \frac{\exp\left(-\frac{1}{\gamma}\sum_{e \in E}\delta_{ep}t_e\right)}{\sum_{\tilde{p} \in P_w}\exp\left(-\frac{1}{\gamma}\sum_{e \in E}\delta_{e\tilde{p}}t_e\right)}, \ p \in P_w.$$

При такой интерпретации связь задачи (5) с логит динамикой, порождающей стохастические равновесия, наиболее наглядна. Последняя формула, в виду того, что $g_p(t) = \sum_{e \in E}\delta_{ep}t_e$ – затраты на пути $p$ на графе $\langle V, E \rangle$, ребра которого взвешены $t$, – есть ни что иное, как отражение следующего принципа поведения (ограниченной рациональности агентов [6]): каждый агент $k$ (пользователь транспортной сети), отвечающий корреспонденции $w \in OD$, выбирает маршрут следования $p \in P_w$, если

$$p = \arg\max_{q \in P_w} \left\{ -g_q(t) + \xi_q^k \right\},$$

где независимые случайные величины $\xi_q^k$, имеют одинаковое двойное экспоненциальное распределение, также называемое распределением Гумбеля[3] [5, 6]:

$$P\left(\xi_q^k < \zeta\right) = \exp\left\{-e^{-\zeta/\gamma - E}\right\}, \ \gamma > 0.$$

Отметим также, что если взять $E \approx 0.5772$ – константа Эйлера, то $M\xi_q^k = 0$, $D\xi_q^k = \gamma^2 \pi^2/6$. Распределение Гиббса получается в пределе, когда число агентов на каждой корреспонденции стремится к бесконечности (случайность исчезает и описание переходит на средние величины). Заметим, что при $\gamma \to 0+$ распределение водителей по путям вырождается, и все водители (агенты) будут использовать только кратчайшие пути

$$-\lim_{\gamma \to 0+}\gamma\psi_w(t/\gamma) = \min_{p \in P_w}g_p(t).$$

Полезно также иметь в виду, что [7]

$$\gamma\psi_w(t/\gamma) = M_{\{\xi_p\}_{p \in P_w}}\left[\max_{p \in P_w}\left\{-g_p(t) + \xi_p\right\}\right].$$

В пределе стабильной динамики задача (5) вырождается в задачу

$$\gamma\psi(t/\gamma) + \langle \bar{f}, t - \bar{t} \rangle \to \min_{t \geq \bar{t}}. \qquad (6)$$

---

[3] Распределение Гумбеля можно объяснить исходя из идемпотентного аналога центральной предельной теоремы (вместо суммы случайных величин – максимум) для независимых случайных величин с экспоненциальным и более быстро убывающим правым хвостом [27]. Распределение Гумбеля возникает в данном контексте, например, если при принятии решения водитель собирает информацию с большого числа разных (независимых) зашумленных источников, ориентируясь на худшие прогнозы по каждому из путей.

В пределе $\gamma \to 0+$ все приведенные выше формулы переходят в соответствующие формулы для модели Бэкмана и модели стабильной динамики [2].

Задачи (5), (6) можно не различать по сложности при композитном подходе к ним [7]. В задаче (6) сепарабельный композит проще сепарабельного композита задачи (5), но зато дополнительно добавляется сепарабельное ограничение простой структуры (для BPR-функций все также сводится к решению уравнения четвертой степени, см. п. 2). В любом случае, основные затраты на каждой итерации при использовании быстрого градиентного метода в композитном варианте [7] связаны с необходимостью расчета градиента гладкой функции[4] $\gamma \psi(t/\gamma)$, при этом константа Липшица градиента этой функции в 2-норме может быть оценена сверху числом [20, 21]

$$\tilde{L}_2 = \frac{1}{\gamma} \sum_{w \in OD} d_w \max_{p \in P_w} \left\| \Theta^{\langle p \rangle} \right\|_2^2.$$

Используя технику быстрого автоматического дифференцирования [28] можно показать, что сложность вычисления градиента этой функции не более чем в четыре раза дороже сложности вычисления значения функции, однако при этом возрастают затраты памяти, поскольку приходится хранить все "дерево вычислений" $\psi(t/\gamma)$.

Приведем, следуя [4, 9], сглаженный идемпотентный аналог метода Форда–Беллмана [29, 30], позволяющий эффективно рассчитывать значение характеристической функции $\psi(t/\gamma)$. Для этого предположим, что любые движения по ребрам графа с учетом их ориентации являются допустимыми, т.е. множество путей, соединяющих заданные две вершины (источник и сток), – это множество всевозможных способов добраться из источника в сток по ребрам имеющегося графа с учетом их ориентации. Этого всегда можно добиться раздутием исходного графа в несколько раз за счет введения дополнительных вершин и ребер. Такое раздутие заведомо можно сделать за $\mathrm{O}(m)$.

Будем считать, что число ребер в любом пути не больше $H = \mathrm{O}(\sqrt{m})$ (диаметр графа с манхетенской структурой, то есть квадратной решетки). Введем классы путей: $P_{ij}^l$ – множество всех путей из $i$ в $j$, состоящих ровно из $l$ ребер, $\tilde{P}_{ij}^l$ – множество всех путей из $i$ в $j$, состоящих из не более чем $l$ ребер. Зафиксируем источник (вершину) $i \in V$, и введем следующие функции для $j \in V$, $l = 1, ..., H$:

$$\begin{cases} a_{ij}^l(t) = \gamma \psi_{P_{ij}^l}(t/\gamma) = \gamma \ln \left( \sum_{p \in P_{ij}^l} \exp\left( -\sum_{e \in E} \delta_{ep} t_e / \gamma \right) \right) \\ b_{ij}^l(t) = \gamma \psi_{\tilde{P}_{ij}^l}(t/\gamma) = \gamma \ln \left( \sum_{p \in \tilde{P}_{ij}^l} \exp\left( -\sum_{e \in E} \delta_{ep} t_e / \gamma \right) \right) \end{cases}.$$

---

[4] Заметим, что возникающие при расчете $\nabla \psi(t/\gamma)$ отношения экспонент стоит сразу же приводить (для большей вычислительной устойчивости) к дробям с числителем равным 1

$$\frac{e^a}{e^a + e^b + ...} = \frac{1}{1 + e^{b-a} + ...}.$$

Некоторые из этих функций могут быть равны $-\infty$. Это означает, что соответствующее множество маршрутов — пустое. Данные функции можно вычислять рекурсивным образом:

$$a_{ij}^1(t) = b_{ij}^1(t) = \begin{cases} -t_e, & e=(i \to j) \in E \\ -\infty, & e=(i \to j) \notin E \end{cases},$$

$$\begin{cases} a_{ij}^{l+1}(t) = \gamma \ln \left( \sum_{k:\, e=(k \to j) \in E} \exp\left( \left( a_{ik}^l(t) - t_e \right)/\gamma \right) \right) \\ b_{ij}^{l+1}(t) = \gamma \ln \left( \exp\left( b_{ij}^l(t)/\gamma \right) + \exp\left( a_{ij}^{l+1}(t)/\gamma \right) \right) \end{cases}, \; j \in V, \; l=1,\ldots,H-1.$$

На каждом шаге $l$ необходимо сделать $\mathrm{O}(m)$ арифметических операций. Следовательно, для вычисления $\psi(t/\gamma)$ необходимо сделать $\mathrm{O}(SHm)$ арифметических операций, где $S = |O|$ — число источников, как правило, можно считать $S \ll m$. Причем вычисление функции $\psi(t/\gamma)$ (и ее градиента) может быть распараллелено на $S$ процессорах. При $\gamma \to 0+$ процедура вырождается в известный метод Форда–Беллмана (динамическое программирование). В процедуре Форда–Беллмана требуется посчитать $H$-степень матрицы $A = \|a_{ij}\|_{i,j \in V}$,

$$a_{ij} = t_e, \; e=(i \to j) \in E;$$
$$a_{ij} = \infty, \; e=(i \to j) \notin E,$$

в идемпотентной математике (вместо обычного поля используется тропическое полуполе [31] со следующими операциями: сложение $a \oplus b = \min\{a,b\}$, произведение $a \otimes b = a+b$). Учитывая, что

$$a \oplus b = \min\{a,b\} = -\lim_{\gamma \to 0+} \gamma \ln \left( \exp(-a/\gamma) + \exp(-b/\gamma) \right),$$
$$a \otimes b = a+b = -\lim_{\gamma \to 0+} \gamma \ln \left( \exp(-a/\gamma) \cdot \exp(-b/\gamma) \right),$$

можно посчитать обычную (над обычным полем) $H$-степень матрицы $A^\gamma = \|a_{ij}^\gamma\|_{i,j \in V}$,

$$a_{ij}^\gamma = e^{-t_e/\gamma}, \; e=(i \to j) \in E;$$
$$a_{ij}^\gamma = 0, \; e=(i \to j) \notin E,$$

и применить поэлементно к полученной матрице $-\gamma \ln(\cdot)$. В пределе $\gamma \to 0+$ получим метод Форда–Беллмана. Однако если не делать предельный переход, то получается нужный нам сглаженный вариант этого алгоритма с такой же временной сложностью. Некоторый аналог этого сглаженного варианта, по сути, и был описан выше.

Оказывается, при определенных условиях (возможности упорядочивания), оценку $\mathrm{O}(SHm)$ можно редуцировать до $\mathrm{O}(Sm)$. Опишем эти условия. Для простоты обозначений будем далее считать, что $S = |D| \ll m$ — число стоков. Переход от источников к стокам также является подготовительным шагом для последующего описания рандомизированного подхода. Предположим, что при любом фиксированном стоке $j \in D$ можно так пронумеровать вершины, что сам сток имеет наибольший номер $v = |V|$ (далее $j = v$) и для любых двух вершин $k, r \in V$ из того, что $k > r$ следует, что не существует ребра из $k$ в $r$. Выполнение этого свойства можно попытаться добиться раздутием исходного графа в несколько раз за счет введения дополнительных вершин и ребер. Такое раздутие и нумерацию можно осуществить за $\mathrm{O}(m)$. Далее, следуя п. 5.4 [32], опишем процедуру расчета $\psi_{iv}(t/\gamma)$, $i \in V$:

$$\psi_{vv}(t/\gamma) \equiv 0,$$

$$\psi_{iv}(t/\gamma) = \gamma \ln\left(\sum_{k:\, e=(i \to k) \in E} \exp\left((\psi_{kv}(t/\gamma) - t_e)/\gamma\right)\right).$$

Отсюда видно, что для вычисления $\psi_{iv}(t/\gamma)$, $i \in V$ необходимо сделать $\mathrm{O}(m)$ арифметических операций, а для вычисления $\psi(t/\gamma)$ соответственно $\mathrm{O}(Sm)$.

Таким образом, используя композитный быстрый градиентный метод [7, 8] можно рассчитывать на следующее число арифметических операций (в действительности, лучше использовать композитный вариант универсального метода [8, 25, 26] поскольку явное использование в методе полученной выше оценки константы Липшица градиента $\tilde{L}_2$ приводит к "лишним итерациям", в виду того, что полученная оценка $\tilde{L}_2$ – завышена; при универсальном подходе метод сам настраивается на подходящую гладкость – адаптивно и не требует никаких констант на входе)

$$\mathrm{O}\left(Sm\sqrt{\frac{\tilde{L}_2 \tilde{R}_2^2}{\varepsilon}}\right) \; \left(\mathrm{O}\left(SHm\sqrt{\frac{\tilde{L}_2 \tilde{R}_2^2}{\varepsilon}}\right) \text{ в общем случае}\right),$$

где $\tilde{R}_2^2$ – евклидов размер решения задачи (5) (или (6) в зависимости от контекста).

**Замечание.** Композитный быстрый градиентный метод (и упоминаемый далее метод зеркального спуска) обладают прямо-двойственной структурой. Это означает, что генерируемая последовательность точек $\{t^k\}$ (в случае композитного быстрого градиентного метода таких последовательности на самом деле три, однако нам понадобятся только две $\{t^k\}$ и $\{\tilde{t}^k\}$ для формулировки результата) обладает следующим свойством

$$\gamma \psi(\tilde{t}^k/\gamma) + \sum_{e \in E} \sigma_e^*(\tilde{t}_e^k) - \min_{t \in \mathrm{dom}\,\sigma^*}\left\{\frac{1}{A_k}\left[\sum_{i=0}^{k} a_i \cdot \left(\gamma \psi(t^i/\gamma) + \langle \nabla \psi(t^i/\gamma), t - t^i \rangle\right)\right] + \sum_{e \in E} \sigma_e^*(t_e)\right\} \leq \frac{C\tilde{L}_2 \breve{R}_2^2}{A_k}, \qquad (7)$$

где константа $C \leq 10$, $a_k \sim k$, $A_k = \sum_{i=0}^{k} a_i$, $A_k \sim k^2$, $\breve{R}_2^2 = \max\{\tilde{R}_2^2, \hat{R}_2^2\}$, $\tilde{R}_2^2 = \frac{1}{2}\|\bar{t} - t^*\|_2^2$ – евклидов размер решения задачи (5), а $\hat{R}_2^2 = \frac{1}{2}\sum_{e \in E}\left(\tau_e(\bar{f}_e^k) - t_e^*\right)^2$ – евклидов размер решения вспомогательной задачи минимизации в (7) (минимум в (7), аналогично задаче (5), достигается во внутренней точки множества $\mathrm{dom}\,\sigma^*$, которое для BPR-функций имеет вид $t \geq \bar{t}$, и лишь при $\mu \to 0+$ минимум может выходить на границу $t_e = \bar{t}_e$), следует сравнить эту оценку с аналогичной оценкой в конце доказательства теоремы 4 [2] и оценкой, приведенной в конце п. 3 [12]. Для метода зеркального спуска правая часть этого неравенства записывается по-другому (см., например, [2]). Положим

$$f^i = -\nabla \psi(t^i/\gamma), \quad x_p^i = d_w \frac{\exp\left(-\frac{1}{\gamma}\sum_{e\in E}\delta_{ep}t_e^i\right)}{\sum_{\tilde{p}\in P_w}\exp\left(-\frac{1}{\gamma}\sum_{e\in E}\delta_{e\tilde{p}}t_e^i\right)}, \quad p \in P_w, \quad \bar{f}^k = \frac{1}{A_k}\sum_{i=0}^{k}a_i f^i, \quad \bar{x}^k = \frac{1}{A_k}\sum_{i=0}^{k}a_i x^i.$$

Учитывая, что

$$\min_{t\in\text{dom}\,\sigma^*}\left\{\frac{1}{A_k}\left[\sum_{i=0}^{k}a_i\langle\nabla\psi(t^i/\gamma),t\rangle\right]+\sum_{e\in E}\sigma_e^*(t_e)\right\} = -\max_{t\in\text{dom}\,\sigma^*}\left\{\left\langle\frac{1}{A_k}\sum_{i=0}^{k}a_i f^i,t\right\rangle - \sum_{e\in E}\sigma_e^*(t_e)\right\} = -\sum_{e\in E}\sigma_e(\bar{f}_e^k),$$

$$\gamma\psi(t^i/\gamma) - \langle\nabla\psi(t^i/\gamma),t^i\rangle = \gamma\psi(t^i/\gamma) + \langle f^i,t^i\rangle = -\sum_{w\in OD}\sum_{p\in P_w}x_p^i\ln(x_p^i/d_w),$$

$$-\sum_{e\in E}\sigma_e(\bar{f}_e^k) - \gamma\frac{1}{A_k}\sum_{i=0}^{k}a_i\sum_{w\in OD}\sum_{p\in P_w}x_p^i\ln(x_p^i/d_w) \leq -\sum_{e\in E}\sigma_e(\bar{f}_e^k) - \gamma\sum_{w\in OD}\sum_{p\in P_w}\bar{x}_p^k\ln(\bar{x}_p^k/d_w),$$

получаем

$$0 \leq \left\{\gamma\psi(\tilde{t}^k/\gamma) + \sum_{e\in E}\sigma_e^*(\tilde{t}_e^k)\right\} + \left\{\sum_{e\in E}\sigma_e(\bar{f}_e^k) + \gamma\sum_{w\in OD}\sum_{p\in P_w}\bar{x}_p^k\ln(\bar{x}_p^k/d_w)\right\} \leq \frac{C\tilde{L}_2\breve{R}_2^2}{A_k},$$

т.е. $\varepsilon$ в приведенных оценках можно понимать как точность по функции решения прямой и двойственной задачи одновременно (можно сказать и точнее: $\varepsilon$ – величина зазора двойственности). Таким образом, критерий останова метода можно завязать на контроль (эффективно вычислимого) зазора двойственности. При этом знание $\breve{R}_2^2$ (и $\tilde{R}_2^2$) методу не требуется (эти параметры не входят в размер шага метода).

Обратим внимание, что эта оценка не зависит от потенциально экспоненциально большого числа путей $|P| = n$ ($n \gg 2^{\sqrt{m}}$, например, для манхетенских сетей), в отличие от оценок, приведенных в пп. 2 – 4. Также обратим внимание, что эта оценка чувствительна к предельному переходу $\gamma \to 0+$, поскольку $\tilde{L}_2 \sim \gamma^{-1}$. В случае малых $\gamma$, можно воспользоваться подходом п. 3 [2], который можно представить для данной задачи следующим образом. Используем рандомизированный метод зеркального спуска, считая вместо $\nabla\psi(t/\gamma)$ стохастический градиент: согласно распределению вероятностей $d_w / \sum_{w\in OD} d_w$ выбираем $w \in OD$ (первый раз это делается за $O(|OD|)$, а все последующие разы за $O(\ln(m))$), затем независимо считаем стохастический градиент $\psi_w(t/\gamma)$. Это можно сделать за $O(m)$ при условиях возможности упорядочивания (см. выше). Опишем, следуя п. 5.4 [32], соответствующую процедуру. Пусть $w = (i,v)$. Считаем за $O(m)$ (см. выше) $\psi_{kv}(t/\gamma)$, $k \in V$. Далее используем следующий подход, базирующийся на формуле полной вероятности. Стартуем из вершины $i$, выбирая ребро, по которому пойдем согласно распределению вероятностей

$$P(e = (i \to k)) = \frac{\exp((\psi_{kv}(t/\gamma) - t_e)/\gamma)}{\exp(\psi_{iv}(t/\gamma)/\gamma)}, \quad k:(i \to k) \in E.$$

По аналогичному принципу выбираем ребро, по которому пойдем из вершины $k$ и т.д. Таким образом, мы сгенерируем случайный путь $p \in P_{iv}$. Поток $\sum_{w \in OD} d_w$ по этому пути определяет потоки на ребрах, входящих в этот путь (на остальных ребрах потоки нулевые). Так задаваемый вектор потоков по ребрам и будет несмещенной оценкой $-\nabla \psi(t/\gamma)$. Оценка сложности (в среднем) метода получается аналогично п. 3 [2]:

$$\tilde{O}\left(m\tilde{M}_2^2 \tilde{R}_2^2 / \varepsilon^2\right),$$

где $\tilde{M}_2^2 = H \cdot \left(\sum_{w \in OD} d_w\right)^2$, а $\tilde{O}(\cdot) = O(\cdot)$ с точностью до логарифмического множителя. Более того, эта оценка имеет и такой же вид, как оценка для аналога этого метода, приведенного в п. 3 [2] при $\gamma = 0$. Описанный подход конкурентно способен с приведенным выше в этом пункте подходом лишь при малых значениях параметра $\gamma > 0$. Можно явно выписывать здесь оценку на $\gamma > 0$ подобно п. 3, но далее мы это сделаем в более интересном случае.

В заключении заметим, что обычные (не стохастические) равновесия можно искать с помощью искусственного введения энтропийной регуляризации. При этом, чтобы с точностью $\varepsilon > 0$ по функции решить исходную задачу, можно действовать следующим образом. Выбрать

$$\gamma = \frac{\varepsilon}{2 \sum_{w \in OD} d_w \ln |P_w|}$$

и решать регуляризованную задачу с точностью $\varepsilon/2$. В этом случае композитный быстрый градиентный метод дает следующую оценку общего числа арифметических операций (в предположении, что транспортный граф имеет правильную упорядоченную структуру, см. выше)

$$\tilde{O}\left(Sm\sqrt{\tilde{M}_2^2 \tilde{R}_2^2 / \varepsilon^2}\right).$$

Мы немного завысили настоящую оценку, чтобы можно было привести получающуюся оценку к максимально удобному для сравнения виду. В зависимости от того, что больше:

$$\tilde{O}\left(m\tilde{M}_2^2 \tilde{R}_2^2 / \varepsilon^2\right) \text{ или } \tilde{O}\left(Sm\sqrt{\tilde{M}_2^2 \tilde{R}_2^2 / \varepsilon^2}\right),$$

следует принимать решение о том вводить или нет регуляризацию.